\documentclass[11pt]{article}
\usepackage{amssymb}
\usepackage{amsmath}
\usepackage{amsbsy}
\newtheorem{lemma}{Lemma}[section]
\newtheorem{proposition}{Proposition}[section]
\newtheorem{theorem}{Theorem}[section]

\newtheorem{remark}{Remark}[section]

\def\proclaim#1{\par \bigskip\noindent {\bf #1}\bgroup\it\ }
\def\endproclaim{\egroup\par\bigskip}

\newbox\TempBox \newbox\TempBoxA
\newcommand{\non}{\nonumber \\}
\def\pr{\textsf{P}} % the symbol P for probability used the sans serif letter
\def\ep{\textsf{E}} % the symbol E for expectation used the sans serif letter
\def\Var{\textsf{Var}} % the symbol Var for covariance used the sans serif letter
\def\Cal#1{{\mathcal #1}}

\def\text#1{\mbox{\rm #1}}
\def\overset#1#2{\stackrel{#1}{#2} }

\def\underwiggle 1{
\ifmmode\setbox\TempBox=\hbox{$ 1$}\else\setbox\TempBox=\hbox{
1}\fi \setbox\TempBoxA=\hbox to \wd\TempBox{\hss\char'176\hss}
\rlap{\copy\TempBox}\smash{\lower9pt\hbox{\copy\TempBoxA}} }

\begin{document}

%\renewcommand{\theequation}
%{\arabic{section}.\arabic{equation}}
%\baselineskip=22pt

%$ $ \vskip 0.8in

%\Large

\title{\huge \bf Precise Asymptotics in Chung's law of the iterated logarithm$^{\ast}$}

%\date{\today}
\date{ }

\maketitle

\begin{center}
{\sc   ZHANG Li-Xin}

\bigskip

{\em Department of Mathematics, Zhejiang University, Hangzhou 310028, China }

\end{center}

\bigskip
\bigskip

%\footnotesize
\noindent $^{\ast}$Research supported by  National Natural Science Foundation of China  (No. 10071072).

\bigskip

\noindent E-mail address: lxzhang {\@} mail.hz.zj.cn

\bigskip
{\rm
%\sf

{\sc Abstract.} \quad Let $X$, $X_1$, $X_2$, $\ldots$ be i.i.d.
random variables with mean zero and positive, finite variance
$\sigma^2$, and set $S_n=X_1+\ldots + X_n$, $n\ge 1$. We prove
that,  if $\ep X^2I\{|X|\ge t\}=o((\log\log t)^{-1})$ as $t\to
\infty$, then for any $a>-1$ and $b>-1$,
\begin{eqnarray*}
 \lim_{\epsilon\nearrow 1/\sqrt{1+a}}
&(\frac{1}{\sqrt{1+a}}-\epsilon)^{b+1} \sum_{n=1}^{\infty}
\frac{(\log n)^a (\log\log n)^b}{n} \pr \Big\{ \max_{k\le
n}|S_k|\le \sqrt{\frac{\sigma^2\pi^2 n}{8\log\log n} }
(\epsilon +a_n)\Big\} \\
&=\frac 4{\pi}(\frac{1}{2(1+a)^{3/2}})^{b+1}\Gamma(b+1),
\end{eqnarray*}
whenever $a_n=o(1/\log\log n)$.

\bigskip



 \vskip 0.2in

\bigskip

}
\newpage

%\newpage

\section{Introduction and main results.}
\setcounter{equation}{0}

Let $\{X, X_n; n\ge 1\}$ be a sequence of i.i.d random variables
with common distribution function $F$, mean $0$ and positive,
finite variance $\sigma^2$, and set $S_n=\sum_{k=1}^nX_k$,
$M_n=\max_{k\le n}|S_k|$, $n\ge 1$. Also let $\log x=\ln(x\vee
e)$, $\log\log x=\log(\log x)$ and $\phi(x)=\sqrt{\pi^2 x
/(8\log\log x)}$.
  Then by the so-called Chung's law of the iterated logarithm
  (LIL) we have
  \begin{equation} \label{eq1.1}
  \liminf_{n\to \infty}\phi(n)M_n=\sigma \quad a.s..
  \end{equation}
This result was first proved by Chung (1948) under
$\ep|X|^3<\infty$, and by Jain and Pruitt (1975) under the sole
assumption of a finite second moment. As pointed by Cs\'aki
(1978), the assumption of a finite second moment is also necessary
for (\ref{eq1.1}) to hold.

As for the usual LIL,  Gut and Sp\u ataru (2000) proved the
following two results on its precise asymptotics.

\proclaim{Theorem A} Suppose that $\ep X=0$, $\ep X^2=\sigma^2$ and $\ep[X^2(\log\log |X|)^{1+\delta}]<\infty$
for some $\delta>0$, and let
$a_n=O(\sqrt{n}/ (\log\log n)^{\gamma})$ for some $\gamma>1/2$. Then
$$\lim_{\epsilon\searrow 1}\sqrt{\epsilon^2-1}\sum_{n=1}^{\infty}\frac 1n
\pr(|S_n|\ge \epsilon\sqrt{2\sigma^2 n\log\log n}+a_n)=1. $$
\endproclaim

\proclaim{Theorem B} Suppose that $\ep X=0$ and $\ep X^2=\sigma^2<\infty$. Then
$$\lim_{\epsilon\searrow 0}\epsilon^2
\sum_{n=1}^{\infty}\frac 1{n\log n}
\pr(|S_n|\ge \epsilon\sqrt{ n\log\log n})=\sigma^2. $$
\endproclaim

The main purpose of this paper is to show similar results on
Chung's LIL under the {\it minimal} conditions by using an
extended Feller's and Einmahl's truncation method.

\begin{theorem} \label{th1}
Let $a>-1$ and $b>-1$ and let $a_n(\epsilon)$ be a function of $\epsilon$ such that
\begin{eqnarray} \label{co0}
a_n(\epsilon) \log \log n \to \tau
\; \text{ as } \;
n\to \infty  \text{ and } \epsilon \nearrow 1/\sqrt{1+a}.
\end{eqnarray}
Suppose that $\ep X=0$, $\ep X^2=\sigma^2<\infty$ and
\begin{eqnarray} \label{co1}
\ep X^2I\{|X|\ge t\}=o((\log\log t)^{-1}) \; \text{ as }\;
t\to \infty.
\end{eqnarray}
Then
\begin{eqnarray}\label{eqT1.1}
& & \lim_{\epsilon\nearrow 1/\sqrt{1+a}}
(\frac{1}{\sqrt{1+a}}-\epsilon)^{b+1} \sum_{n=1}^{\infty}
\frac{(\log n)^a (\log\log n)^b}{n}
\pr \Big\{M_n \le \sigma \phi(n) (\epsilon
+a_n(\epsilon))\Big\}
\non
& & \qquad \qquad \qquad=\frac
4{\pi}(\frac{1}{2(1+a)^{3/2}})^{b+1}\Gamma(b+1)\exp\{2(1+a)^{3/2}\tau\}.
\end{eqnarray}
Here, $\Gamma(\cdot)$ is a gamma function.
Conversely, if (\ref{eqT1.1}) holds for some $a>-1$, $b>-1$ and $0<\sigma<\infty$, then $\ep X=0$, $\ep X^2=\sigma^2$ and
\begin{equation}\label{eqT1.2}
\liminf_{t\to \infty}(\log\log t)\ep X^2I\{|X|\ge t\}=0.
\end{equation}
\end{theorem}

\begin{theorem} \label{th2}
Suppose that $\ep X=0$ and $\ep X^2=\sigma^2<\infty$.
For $b>-1$, we have
\begin{eqnarray}\label{eqT2.1}
& & \lim_{\epsilon\nearrow \infty}
\epsilon^{-2(b+1)} \sum_{n=1}^{\infty}
\frac{(\log\log n)^b}{n\log n} \pr \Big\{M_n \le
\sigma\phi(n) \epsilon
\Big\}
\non
& & \qquad \qquad \qquad=\frac
4{\pi}\Gamma(b+1)\sum_{k=0}^{\infty}\frac{(-1)^k}{(2k+1)^{2b+3}}.
\end{eqnarray}
Conversely, if (\ref{eqT2.1}) holds for some $b>-1$ and $0<\sigma<\infty$, then $\ep X=0$ and  $\ep X^2=\sigma^2$.
\end{theorem}

Theorems \ref{th1} and \ref{th2} are also related to the integral
test which refines (\ref{eq1.1}).
The first one of the results on the integral test for (\ref{eq1.1})
is due to Chung (1948) who obtained that
if $\ep X=0$,  $\ep X^2=\sigma^2$ and $\ep |X|^3<\infty$,
then for any eventually non-increasing
$\psi:[1,\infty)\to (0, \infty)$,
\begin{eqnarray} \label{eq1.2}
& & \pr\big( M_n\le  \sqrt{\sigma^2 \pi^2 n/8}\psi(n) \;
i.o.\big)=0 \text{ or } =1 \non
 & & \text{ according as }
J(\psi):=\sum_{n=1}^{\infty}\frac{1}{n\psi(n)^2}\exp(-1/\psi(n)^2)<\infty
\text{ or } = \infty.
\end{eqnarray}
Einmahl (1993) proved (\ref{eq1.2}) under the  minimal condition
that
\begin{eqnarray} \label{co2}
\ep X^2I\{|X|\ge t\}=O((\log\log t)^{-1}) \; \text{ as }\; t\to
\infty.
\end{eqnarray}

Our next theorem gives a result on a convergence rate of
(\ref{eq1.1}) and (\ref{eq1.2}).
\begin{theorem} \label{th3}
Let $a>-1$ and $b>-1$.  Suppose that $\ep X=0$,
 $\ep X^2=\sigma^2$, and that the condition (\ref{co2}) is
satisfied, then for any eventually non-increasing
$\psi:[1,\infty)\to (0, \infty)$,
\begin{eqnarray} \label{eqT3.1}
& & \sum_{n=1}^{\infty} \frac{(\log n)^a (\log\log
n)^b}{n}\pr\big( M_n\le  \sqrt{\sigma^2 \pi^2 n/8} \psi(n)
\big)<\infty \text{ or } =\infty \non
 & & \text{ according as }
J_{ab}(\psi):=\sum_{n=1}^{\infty}\frac{(\log n)^a (\log\log
n)^b}{n}\exp(-1/\psi(n)^2)<\infty \text{ or } = \infty.
\end{eqnarray}
Conversely, if (\ref{eqT3.1}) holds for some $a>-1$, $b>-1$,  $0<\sigma<\infty$ and any eventually non-increasing
$\psi(x)$,
then $\ep X=0$, $\ep X^2=\sigma^2$ and
\begin{equation}\label{eqT3.2}
\liminf_{t\to \infty}(\log\log t)\ep X^2I\{|X|\ge t\}<\infty.
\end{equation}
\end{theorem}

By (\ref{eqT3.1}), we know that the infinite series in (\ref{eqT1.1})
converges whenever $\epsilon<1/\sqrt{1+a}$, and diverges whenever
$\epsilon>1/\sqrt{1+a}$.

\begin{remark}
Note that the conditions (\ref{co1}) and (\ref{co2}) is
sharp. A sufficient condition for them is given by
$$ \ep X^2 \log\log |X|<\infty. $$
However, when $a>0$, the sufficient and necessary condition for
$$\sum_{n=1}^{\infty} \frac{(\log n)^a (\log\log
n)^b}{n}\pr\big( M_n> \epsilon \sqrt{2\sigma^2n\log\log n}
\big)<\infty, \quad \epsilon>\sqrt{1+a}, $$
to hold is that
$$ \ep  X^2 (\log|X|)^a(\log\log |X|)^{b-1}<\infty. $$
\end{remark}

\begin{remark} By using a result of Einmahl (1987) instead of our
Lemma \ref{lem5}, one can extend Theorems \ref{th1}-\ref{th3} to
multidimensional random variables.
\end{remark}

The proofs of Theorems \ref{th1}-\ref{th3} are given in Section 4.
Before that, we first verify (\ref{eqT1.1}), (\ref{eqT2.1}) and
(\ref{eqT3.1}) under the assumption that $F$ is the normal
distribution in Section 2, after which, by using the truncation
and approximation method, we then show that the probabilities in
(\ref{eqT1.1}), (\ref{eqT2.1}) and (\ref{eqT3.1}) can be replaced
by those for normal random variables in Section 3. Throughout this
paper, we let $K(\alpha,\beta,\cdots)$, $C(\alpha,\beta,\cdots)$
etc denote positive constants which depend on $\alpha,\beta,
\cdots$ only, whose values can differ in different places.
$a_n\sim b_n$ means that $a_n/b_n\to 1$.

%%%%%%%%%%%%%%%%%%%%%%%%%%%%%%%%%%%%%%
%%%%%%%%%%%%%%%%%%%%%%%%%%%%%%%%%%%%%%%%

\section{Normal cases.}
\setcounter{equation}{0}

In this section, we prove Theorems \ref{th1}-\ref{th3} in the case that $\{X, X_n; n\ge 1\}$ are normal random variables.
Let $\{W(t); t\ge 0\}$ be a standard Wiener process. Our results are as follows.

\begin{proposition}\label{prop2.1}
Let $a>-1$ and $b>-1$ and let $a_n(\epsilon)$ be a function of $\epsilon$ satisfying (\ref{co0}).
Then
\begin{eqnarray}\label{eqprop2.1.1}
&& \lim_{\epsilon\nearrow 1/\sqrt{1+a}}
(\frac{1}{\sqrt{1+a}}-\epsilon)^{b+1} \sum_{n=1}^{\infty}
\frac{(\log n)^a (\log\log n)^b}{n}
 \non & & \qquad \qquad \qquad
 \qquad \qquad \cdot \pr \Big\{\sup_{0\le s\le 1}|W(s)| \le
\sqrt{\frac{\pi^2 }{8\log\log n} } (\epsilon +a_n(\epsilon))\Big\}
\non & & \qquad \qquad \qquad=\frac
4{\pi}(\frac{1}{2(1+a)^{3/2}})^{b+1}\Gamma(b+1)
\exp\Big\{2(1+a)^{3/2}\tau\Big\}.
\end{eqnarray}
\end{proposition}

\begin{proposition}\label{prop2.2}
For any $b>-1$, we have
\begin{eqnarray*}
& & \lim_{\epsilon\nearrow \infty} \epsilon^{-2(b+1)}
\sum_{n=1}^{\infty} \frac{(\log\log n)^b}{n\log n}
 \pr\Big\{\sup_{0\le s\le 1}|W(s)| \le \epsilon \sqrt{\frac{\pi^2 }{8\log\log
n} } \Big\}
\\
& & \qquad \qquad \qquad
=\frac4{\pi}\Gamma(b+1)
\sum_{k=0}^{\infty}\frac{(-1)^k}{(2k+1)^{2b+3}}.
\end{eqnarray*}
\end{proposition}

\begin{proposition} \label{prop2.3}
For any $a>-1$ and $b>-1$, we have that
\begin{eqnarray*}
& & \sum_{n=1}^{\infty} \frac{(\log n)^a (\log\log
n)^b}{n}\pr\big( \sup_{0\le s\le 1}|W(s)|\le \sqrt{\pi^2/8} \psi(n) \big)<\infty
\text{ or } =\infty
\\
 & & \text{ according as }
J_{ab}(\psi)<\infty \text{ or } = \infty.
\end{eqnarray*}
\end{proposition}

The following lemma will be used in the proofs.

\begin{lemma}\label{lem2.1}
Let $\{W(t); t\ge 0\}$ be a standard Wiener process.
Then for all $x>0$,
\begin{equation}\label{eqL2.1.1}
 \pr \big(\sup_{0\le s\le 1}|W(s)| \le x\big)=
\frac4{\pi}\sum_{k=1}^{\infty}
\frac{(-1)^k}{2k+1}\exp\Big\{-\frac{\pi^2(2k+1)^2}{8x^2}\Big\}
\end{equation}
and
$$
\pr\big( \sup_{0\le s\le 1}|W(s)|\le x \big)
\sim \frac 4{\pi}\exp\big\{-\frac{\pi^2}{8 x^2} \big\}
\; \text{ as } \; x\to 0.
$$
\end{lemma}
{\bf Proof.} It is well known. See Ciesielski and Taylor (1962).

\bigskip

Now, we turn to prove the propositions.

{\noindent\bf Proof Proposition \ref{prop2.1}:} First, note that
the limit in (\ref{eqprop2.1.1}) does not depend on any finite
terms of the infinite series.
Secondly, by Lemma \ref{lem2.1} and the condition ({\ref{co0}) we
have
\begin{eqnarray*}
& &\pr \Big\{\sup_{0\le s\le 1}|W(s)| \le
\sqrt{\frac{\pi^2 }{8\log\log n} }
(\epsilon +a_n(\epsilon))\Big\}
\sim \frac 4{\pi}\exp\Big\{-\frac{\log\log n}{(\epsilon+a_n(\epsilon))^2}\Big\} \\
& & \qquad= \frac 4{\pi}\exp\Big\{-\frac{\log\log
n}{\epsilon^2+2\epsilon a_n(\epsilon)+a_n^2(\epsilon)}\Big\}\\
& &\qquad \sim \frac 4{\pi}\exp\Big\{-\frac{\log\log
n}{\epsilon^2}\Big\}
\exp\Big\{\frac{2}{\epsilon^3}a_n(\epsilon)\log\log n \Big\}
\end{eqnarray*}
as $n\to \infty$, uniformly in
$\epsilon\in (1/\sqrt{1+a}-\delta,1/\sqrt{1+a})$
for some $\delta>0$.
So, for any $0<\theta<1$,
there exist $\delta>0$ and $n_0$ such that for all $n\ge n_0$ and
$\epsilon\in (1/\sqrt{1+a}-\delta,1/\sqrt{1+a})$,
$$\begin{array}{ll}
& \frac 4{\pi}\exp\Big\{-\frac{\log\log n}{\epsilon^2}\Big\}
\exp\Big\{2(1+a)^{3/2}\tau-\theta\Big\}\\
\le &\pr \Big\{\sup_{0\le s\le 1}|W(s)| \le
  \sqrt{\frac{\pi^2 }{8\log\log n} }
  ( \epsilon +a_n(\epsilon))\Big\} \\
\le & \frac 4{\pi}\exp\Big\{-\frac{\log\log n}{\epsilon^2}\Big\}
\exp\Big\{2(1+a)^{3/2}\tau+\theta\Big\},
\end{array} $$
by the condition ({\ref{co0}) again. Also,
\begin{eqnarray*}
& &  \lim_{\epsilon\nearrow 1/\sqrt{1+a}}
(\frac{1}{\sqrt{1+a}}-\epsilon)^{b+1} \sum_{n=1}^{\infty}
\frac{(\log n)^a (\log\log n)^b}{n}
\exp\Big\{-\frac{\log\log n}{\epsilon^2}\Big\} \\
&=& \lim_{\epsilon\nearrow 1/\sqrt{1+a}}
(\frac{1}{\sqrt{1+a}}-\epsilon)^{b+1}
\int_{e^e}^{\infty} \frac{(\log x)^a (\log\log x)^b}{x}
\exp\Big\{-\frac{\log\log x}{\epsilon^2}\Big\} dx \\
&=&
 \lim_{\epsilon\nearrow 1/\sqrt{1+a}}
(\frac{1}{\sqrt{1+a}}-\epsilon)^{b+1}
\int_1^{\infty} y^b
\exp\Big\{-y\big(\frac{1}{\epsilon^2}-1-a\big)\Big\} dy \\
&=& \lim_{\epsilon\nearrow 1/\sqrt{1+a}}
(\frac{1}{\sqrt{1+a}}-\epsilon)^{b+1}
\big(\frac{1}{\epsilon^2}-1-a\big)^{-b-1}
\int_{\frac{1}{\epsilon^2}-1-a}^{\infty} y^b
e^{-y} dy \\
&=& (\frac{1}{2(1+a)^{3/2}})^{b+1} \int_0^{\infty}
y^be^{-y} dy
 =(\frac{1}{2(1+a)^{3/2}})^{b+1}\Gamma(b+1).
\end{eqnarray*}
(\ref{eqprop2.1.1}) is now proved.

\bigskip

{\noindent\bf Proof Proposition \ref{prop2.2}:} Notice that (\ref{eqL2.1.1}),
and for any $m\ge 1$ and all $x>0$,
\begin{eqnarray*}
& &\frac4{\pi}\sum_{k=0}^{2m+1}
\frac{(-1)^k}{2k+1}\exp\Big\{-\frac{\pi^2(2k+1)^2}{8x^2}\Big\}
\\
&\le& \pr \big(\sup_{0\le s\le 1}|W(s)| \le x\big)\\
&\le&
\frac4{\pi}\sum_{k=0}^{2m}
\frac{(-1)^k}{2k+1}\exp\Big\{-\frac{\pi^2(2k+1)^2}{8x^2}\Big\}.
\end{eqnarray*}
It is sufficient to show that for any $q>0$,
\begin{eqnarray}\label{eq2.3}
 \lim_{\epsilon\nearrow \infty}
\epsilon^{-2(b+1)} \sum_{n=1}^{\infty}
\frac{(\log\log n)^b}{n\log n}
\exp\Big\{-q\frac{\log\log n}{\epsilon^2}\Big\}
=\Gamma(b+1)q^{-(b+1)}.
\end{eqnarray}
Now, the left hand of the above equality equals
\begin{eqnarray*}
 & &\lim_{\epsilon\nearrow \infty}
\epsilon^{-2(b+1)} \int_{e^e}^{\infty}
\frac{(\log\log x)^b}{x\log x}
\exp\Big\{-q\frac{\log\log x}{\epsilon^2}\Big\}dx
\\
&=&\lim_{\epsilon\nearrow \infty}
\epsilon^{-2(b+1)} \int_1^{\infty}y^b
\exp\Big\{-y\frac{q}{\epsilon^2}\Big\}dy \\
&=&\lim_{\epsilon\nearrow \infty}
q^{-(b+1)} \int_{q/\epsilon^2}^{\infty}y^be^{-y}dy
=\Gamma(b+1)q^{-(b+1)}.
\end{eqnarray*}
The proposition is proved.

\bigskip
{\noindent\bf Proof Proposition \ref{prop2.3}:} It is obvious, since
$$ \frac 2{\pi}\exp\Big\{-\frac{\pi^2}{ 8 x^2}\Big\}
\le \pr\{\sup_{0\le s\le 1}|W(s)|\le x\}\le
\frac 4{\pi}\exp\Big\{-\frac{\pi^2}{ 8 x^2}\Big\}, \quad
\forall x>0. $$

%%%%%%%%%%%%%%%%%%%%%%%%%%%%%%%%%%%%%
%%%%%%%%%%%%%%%%%%%%%%%%%%%%%%%%%%%%

\section{Truncation and Approximation.}
\setcounter{equation}{0}

The purpose of this section is to extend Feller's and Einmahl's
truncation methods and to show that the probabilities in
(\ref{eqT1.1}), (\ref{eqT2.1}) and (\ref{eqT3.1}) for $M_n$ can be
approximated by those for $\sqrt{n}\sup_{0\le s\le 1}|W(s)|$.

Suppose that $\ep X=0$ and $\ep X^2=\sigma^2<\infty$.
Without losing of generality, we assume that $\sigma=1$.
 Let $0<p<1/2$. For each $n$ and $1\le j\le n$,
we let $X_{nj}^{\prime}=X_{nj}I\{|X_j|\le \sqrt{n}/\log^p n \}$,
$X_{nj}^{\ast}=X_{nj}^{\prime}-\ep [X_{nj}^{\prime}]$,
$S_{nj}^{\prime}=\sum_{i=1}^jX_{nj}^{\prime}$,
$S_{nj}^{\ast}=\sum_{i=1}^jX_{nj}^{\ast}$,
$M_n^{\ast}=\max_{k\le n}|S_{nk}^{\ast}|$ and
$B_n=\sum_{k=1}^n \Var(X_{nk}^{\ast})$.
The following two propositions are the main results of this section.

\begin{proposition}\label{prop3.1}
Let $a>-1$, $b>-1$ and $0<p<1/2$. Then there exist $\delta>0$ and a sequence of positive numbers $\{q_n \}$ such that
\begin{eqnarray} \label{eqprop3.1.1}
&&\pr\Big\{\sup_{0\le s\le 1}|W(s)| \le \epsilon
\sqrt{\frac{\pi^2
}{8\log\log n} }
 - \frac{3}{(\log\log n)^2} \Big\}
-q_n \non &\le& \pr\Big\{M_n \le \epsilon \sqrt{\frac{\pi^2
B_n}{8\log\log n} } \Big\} \non &\le&\pr\Big\{\sup_{0\le s\le
1}|W(s)| \le \epsilon \sqrt{\frac{\pi^2 }{8\log\log n} }
 +
\frac{3}{(\log\log n)^2} \Big\} +q_n, \non & &  \qquad \forall
\epsilon\in (\frac 1{\sqrt{1+a}}-\delta, \frac
1{\sqrt{1+a}}+\delta), \quad n\ge 1
\end{eqnarray}
and
\begin{eqnarray}\label{eqprop3.1.2}
\sum_{n=1}^{\infty}\frac{(\log n)^a(\log\log n)^b}{n}
q_n \le K(a,b,p,\delta)<\infty.
\end{eqnarray}
\end{proposition}

\begin{proposition}\label{prop3.2}
If $n$ is large enough, then
\begin{eqnarray}\label{eqprop3.2.1}
& &
\pr\big(\sup_{0\le s\le 1}|W(s)|\le x-3/(\log\log n)^2\big)-q_n^{\ast}
\le \pr\big(M_n\le x\sqrt{B_n}\big)
\non
&\le&
\pr\big(\sup_{0\le s\le 1}|W(s)|\le x+3/(\log\log n)^2\big)+q_n^{\ast},
\quad
\forall x>0,
\end{eqnarray}
where $q_n^{\ast}\ge 0$ satisfies
\begin{equation} \label{eqprop3.2.2}
\sum_{n=1}^{\infty} \frac{(\log\log n)^b}{n\log n} q_n^{\ast} \le
K(b,p)<\infty.
\end{equation}
\end{proposition}

To show this two results, we need some lemmas.

\begin{lemma}\label{lem1}
For any $x>0$ and $0<\delta<1$,
there exists a positive constant $C=C(x,\delta)$ such that
\begin{itemize}
\item[{(a)}] $C^{-1}\exp\{-\frac{1+\delta}{x^2}\log\log n\}
 \le \pr\{M_n^{\ast}\le x\phi(n)\}
 \le C\exp\{-\frac{1-\delta}{x^2}\log\log n\}$,
\item[{(b)}] $C^{-1}\exp\{-\frac{1+\delta}{x^2}\log\log n\}
 \le \pr\{M_n \le x\phi(n)\}
 \le C\exp\{-\frac{1-\delta}{x^2}\log\log n\}$.
\end{itemize}
\end{lemma}

{\bf Proof.}  This lemma is so-called small deviation theorem. It
follows from Theorem 2 of Shao (1995) by noting that $B_n\sim n$.
(See also Shao 1991).

\begin{lemma}\label{lem2}
For any $x>0$, $A>0$ and $0<\delta<1$, there exists a  positive constant $C=C(x,\delta)$ such that
$$ \pr\{\max_{k\le q}|S_k+z_1|\vee
\max_{q<k\le n}|S_k+z_2|\le x\phi(n)\} \le
C\exp\big\{-\frac{1-\delta}{x^2}\log\log n\big\} $$ holds
uniformly in $|z_1|\le A\phi(n)$, $|z_2|\le A\phi(n)$ and $1\le
q\le n$.
\end{lemma}

{\bf Proof.} Without losing of generality,
we can assume that $0<\delta<\frac 1{8^8 5!}$.
We follow the lines of the proof of (17) in Shao (1995).
Write $x_n=x\phi(n)$. Put $M=\delta^{-2}$.
For fixed $n$, define $m_0=0$,
$$ m_i=\max\{j: j\le Mi x_n^2\},
 \quad \text{ for } i\le l:=\max\{i: m_i\le n-1\} $$
and $m_{l+1}=n$. It is easily seen that
$$(1-\delta/4)M x_n^2\le m_i-m_{i-1}\le (1+\delta/4) M x_n^2 $$
and
$$\frac{n}{M x_n^2}-1\le l \le \frac{n}{M x_n^2} $$
provided $n$ is sufficiently large.
From Lemmas 3 and 1 of Shao (1995) and the Anderson's inequality,
it follows that, there exists an integer $n_0$ such that
$\forall n\ge n_0$, $\forall 1\le j\le l$, $\forall |y|\le x_n+A\phi(n)$,
$\forall |y_j|\le A\phi(n)$,
\begin{eqnarray*}
& & \pr\Big(\max_{k\le m_j-m_{j-1}}|S_{m_{j-1}+k}-S_{m_{j-1}}+y_j+y|\le x_n\Big) \\
&\le & e^{-3M}
 +\pr\Big(\sup_{0\le s\le 1}|W(s)+(y_j+y)(m_j-m_{j-1})^{-1/2}|\le x_n(m_j-m_{j-1})^{-1/2}\Big) \\
&\le & e^{-3M}
 +\pr\Big(\sup_{0\le s\le 1}|W(s)|\le x_n(m_j-m_{j-1})^{-1/2}\Big) \\
&\le & e^{-3M}
+4\exp\Big(-\frac{\pi^2(m_j-m_{j-1})}{8x_n^2}\Big)
\le e^{-3M}
+4\exp\Big(-\frac{\pi^2 M(1-\delta/4)}{8}\Big) \\
&\le &\frac 12e^{-2M}
+\frac 12 \exp\Big(-\frac{\pi^2 M(1-\delta/2)}{8}\Big)
\le \exp\Big(-\frac{\pi^2 M(1-\delta/2)}{8}\Big),
\end{eqnarray*}
where $\{W(t); t\ge 0\}$ is a standard Wiener process. Obviously,
there exists an $i$ such that $m_{i-1}<q\le m_i$. Let $y_j=z_1$ if
$j\le i-1$, and $z_2$ if $j\ge i$. Then
\begin{eqnarray*}
& & \pr\big( \max_{k\le q}|S_k+z_1|\vee
  \max_{q<k\le n}|S_k+z_2|\le x\phi(n)\big)
\le \pr\big( \max_{j\le l,j\ne i}\max_{m_{j-1}<k\le m_j}|S_k+y_j|
 \le x_n \big) \\
&= &\ep\Big\{I\{ \max_{j\le l-1,j\ne i}\max_{m_{j-1}<k\le m_j}|S_k+y_j|
 \le x_n\}I\{ \max_{m_{l-1}<k\le m_l, l\ne i}|S_k+y_j|
 \le x_n\} \Big\} \\
&\le &\ep\Big\{I\{ \max_{j\le l-1,j\ne i}\max_{m_{j-1}<k\le m_j}|S_k+y_j|
 \le x_n\}\\
& & \qquad \times \ep\big[I\{ \max_{m_{l-1}<k\le m_l}|S_k+y_j|
 \le x_n\}|S_k, k\le m_{l-1}\big] \Big\} \\
&= &\int_{-x_n-y_{l-1}}^{x_n-y_{l-1}}
 \pr\Big( \max_{m_{l-1}<k\le m_l}|S_k-S_{m_{l-1}}+y+y_l|
 \le x_n\Big) \\
& & \qquad \quad d \pr\Big( \max_{j\le l-1,j\ne i}\max_{m_{j-1}<k\le m_j}|S_k+y_j| \le x_n, S_{m_{l-1}}<y\Big) \\
&\le& \exp\Big(-\frac{\pi^2 M(1-\delta/2)}{8}\Big)
\pr\Big( \max_{j\le l-1,j\ne i}\max_{m_{j-1}<k\le m_j}|S_k+y_j|
   \le x_n\Big) \\
&\le &\ldots \le \exp\Big(-\frac{\pi^2 M(1-\delta/2)(l-1)}{8}\Big) \\
&\le & C \exp\Big(-\frac{(1-\delta)\pi^2 n}{8x_n^2}\Big)
\le C \exp\Big(-\frac{(1-\delta)}{x^2}\log\log n\Big).
\end{eqnarray*}
For $n\le n_0$, it is obvious that
$$\pr\big( \max_{k\le q}|S_k+z_1|\vee
  \max_{k<k\le n}|S_k+z_2|\le x\phi(n)\big) \le 1
\le C \exp\Big(-\frac{(1-\delta)}{x^2}\log\log n\Big). $$
Lemma 2 is proved.

\begin{lemma}\label{lem3}
Define $\Delta_n=\max_{k\le n}|S_{nk}^{\ast}-S_k|$.
Let $a>-1$, $b>-1$ and $0<p<1/2$.
Then there exist constants $\delta_0=\delta_0(a,p)>0$
and $K=K(a,b,p)$ such that $\forall 0<\delta<\delta_0$,
$$\sum_{n=1}^{\infty} \frac{(\log n)^a(\log\log n)^b}{n}I_n\le K \ep X^2<\infty, $$
where
$$I_n=\pr\Big(\Delta_n\ge \sqrt{n}/(\log\log n)^2,
 M_n^{\ast}\le \phi(n)\big(\frac 1{\sqrt{1+a}}+\delta\big)\Big). $$
\end{lemma}
{\bf Proof.} It is sufficient to show that
\begin{eqnarray}\label{eqL3.1}
\sum_{n=1}^{\infty} \frac{(\log n)^a(\log\log n)^b}{n}
\pr\Big(\Delta_n\ge \sqrt{n}/(\log\log n)^2,
 M_n^{\ast}\le \frac{\phi(n)}{\sqrt{1+a-\delta}}\Big)
\le C \ep X^2,
\end{eqnarray}
whenever $0<\delta<1+a$ and $0<\delta<1-2p$.
Let $\beta_n=n \ep[|X|I\{|X|>\sqrt{n}/\log^p n\}]$. Then
$|\ep\sum_{i=1}^jX_{ni}^{\prime}|\le \beta_n$, $1\le j\le n$.
Setting
$$\Cal L=\{n:\beta_n\le \frac 18 \sqrt{n}/(\log\log n)^2\},$$
we have
$$\{\Delta_n\ge \sqrt{n}/(\log\log n)^2\}
\subset
\bigcup_{j=1}^n \{ X_j\ne X_{nj}^{\prime}\}, \quad n\in \Cal L. $$
So for $n\in \Cal L$,
\begin{eqnarray*}
I_n^{\prime}:
&=&\pr\Big(\Delta_n\ge \sqrt{n}/(\log\log n)^2,
 M_n^{\ast}\le \frac{\phi(n)}{\sqrt{1+a-\delta}}\Big)\\
&\le& \sum_{j=1}^n \pr\Big(X_j\ne X_{nj}^{\prime}, M_n^{\ast}\le \frac{\phi(n)}{\sqrt{1+a-\delta}}\Big).
\end{eqnarray*}
Observer that $X_{nj}^{\prime}=0$ whenever $X_j\ne X_{nj}^{\prime}$,
$j\le n$, so that by Lemma \ref{lem1}(a) we have for $n$ large enough and all $1\le j\le n$,
\begin{eqnarray*}
& &\pr\Big(X_j\ne X_{nj}^{\prime},
M_n^{\ast}\le \frac{\phi(n)}{\sqrt{1+a-\delta}}\Big)\\
&=&\pr\Big(X_j\ne X_{nj}^{\prime},
\max_{k\le j-1}|S_{nk}^{\ast}|\vee
\max_{j<k\le n}|S_{nk}^{\ast}-X_{nj}^{\prime}|
\le \frac{\phi(n)}{\sqrt{1+a-\delta}}\Big)\\
&=&\pr(X_j\ne X_{nj}^{\prime})
\pr\Big(\max_{k\le j-1}|S_{nk}^{\ast}|\vee
\max_{j<k\le n}|S_{nk}^{\ast}-X_{nj}^{\prime}|
\le \frac{\phi(n)}{\sqrt{1+a-\delta}}\Big)\\
&\le &\pr\Big(X_j\ne X_{nj}^{\prime})
\pr\Big(M_n^{\ast}\le \frac{\phi(n)}{\sqrt{1+a-\delta}}
+|X_{nj}^{\prime}|\Big)\\
&\le &\pr(|X|>\sqrt{n}/\log^p n)
\pr\Big(M_n^{\ast}\le \frac{\phi(n)}{\sqrt{1+a-\delta}}
+\sqrt{n}/\log^p n\Big)\\
&\le &C\pr(|X|>\sqrt{n}/\log^p n)\exp\{(-1-a+\delta+\delta^{\prime})\log\log n\},
\end{eqnarray*}
where $0<\delta^{\prime}<1-2p-\delta$. So,
\begin{eqnarray*}
& & \sum_{n\in \Cal L} \frac{(\log n)^a (\log\log n)^b }{n}
I_n^{\prime}
\le C\sum_{n=1}^{\infty}\pr(|X|>\sqrt{n}/\log^p n)(\log n)^{ \delta+\delta^{\prime}-1 }(\log\log n)^b \\
&\le&  \sum_{n=1}^{\infty}\sum_{j=n}^{\infty}
 \pr\big(\sqrt{j}/\log^p j<|X|\le \sqrt{j+1}/\log^p (j+1)\big)
(\log n)^{\delta+\delta^{\prime}-1}(\log\log n)^b \\
&\le&  \sum_{j=1}^{\infty}
 \pr\big(\sqrt{j}/\log^p j<|X|\le \sqrt{j+1}/\log^p (j+1)\big)
\sum_{n=1}^j (\log n)^{\delta+\delta^{\prime}-1 }(\log\log n)^b \\
&\le&  \sum_{j=1}^{\infty}
 \pr\big(\sqrt{j}/\log^p j<|X|\le \sqrt{j+1}/\log^p (j+1)\big)
j (\log j)^{\delta+\delta^{\prime}-1 }(\log\log j)^b \\
&\le& C\ep\Big[X^2(\log |X|)^{\delta+\delta^{\prime}+2p-1 }(\log\log |X|)^b \Big]
\le C\ep X^2.
\end{eqnarray*}
If $n\not\in \Cal L$, then by Lemma \ref{lem1}(a) we have
\begin{eqnarray*}
 I_n^{\prime}
&\le& \pr\Big(M_n^{\ast}\le  \frac{\phi(n)}{\sqrt{1+a-\delta}}\Big)
\le C\exp\{(-1-a+\delta+\delta^{\prime})\log\log n\}\\
&=& C(\log n)^{-1-a+\delta+\delta^{\prime}}.
\end{eqnarray*}
It follows that
\begin{eqnarray*}
& & \sum_{n\not\in \Cal L}\frac{(\log n)^a(\log\log
n)^b}{n}I_n^{\prime} \le C\sum_{n\not\in \Cal L}
\frac{(\log n)^{\delta+\delta^{\prime}-1}(\log\log n)^b}{n} \\
&\le& 8C\sum_{n\not\in \Cal L}
\frac{(\log n)^{\delta+\delta^{\prime}-1}(\log\log n)^{b+2}}{n^{3/2}}
 \beta_n  \\
&\le& 8C\sum_{n=1}^{\infty}\frac{(\log n)^{\delta+\delta^{\prime}-1}(\log\log n)^{b+2}}{n^{1/2}}
\\
& & \qquad \cdot \sum_{j=n}^{\infty}
\ep\Big[|X|I\{\sqrt{j}/\log^p j<|X|\le \sqrt{j+1}/\log^p(j+1)\}\Big]\\
&=& 8C\sum_{j=1}^{\infty} \ep\Big[|X|I\{\sqrt{j}/\log^p j<|X|\le
\sqrt{j+1}/\log^p(j+1)\}\Big] \\
& & \qquad \cdot \sum_{n=1}^j
\frac{(\log n)^{\delta+\delta^{\prime}-1}(\log\log n)^{b+2}}{n^{1/2}}\\
&\le& C\sum_{j=1}^{\infty} \ep\Big[|X|I\{\sqrt{j}/\log^p j<|X|\le
\sqrt{j+1}/\log^p(j+1)\}\Big]\\
& & \qquad \cdot
\sqrt{j}(\log j)^{\delta+\delta^{\prime}-1}(\log\log j)^{b+2}\\
&\le& C\ep\Big[X^2(\log |X|)^{\delta+\delta^{\prime}-1+p}
(\log\log |X|)^{b+2}\Big]\le C \ep X^2.
\end{eqnarray*}
(\ref{eqL3.1}) is proved.

\begin{lemma}\label{lem4}
Let $a>-1$, $b>-1$ and $0<p<1/2$.
Then there exist constants $\delta_0=\delta_0(a,p)>0$
and $K=K(a,b,p)$ such that $\forall 0<\delta<\delta_0$,
$$\sum_{n=1}^{\infty} \frac{\log^q n(\log\log n)^b}{n}II_n\le K \ep X^2<\infty, $$
where
$$II_n=\pr\Big(\Delta_n\ge \sqrt{n}/(\log\log n)^2,
 M_n\le \phi(n)\big(\frac 1{\sqrt{1+a}}+\delta\big)\Big). $$
\end{lemma}
{\bf Proof.} It is enough to show that
\begin{eqnarray*}
\sum_{n=1}^{\infty} \frac{\log^q n(\log\log n)^b}{n}
\pr\Big(\Delta_n\ge \sqrt{n}/(\log\log n)^2,
 M_n\le \frac{\phi(n)}{\sqrt{1+a-\delta}}\Big)
\le C \ep X^2,
\end{eqnarray*}
whenever $0<\delta<1+a$ and $0<\delta<1-2p$.
Let $\beta_n$ and $\Cal L$ be defined as in the proof of Lemma \ref{lem3}.
Then for $n\in \Cal L$,
\begin{eqnarray*}
\pr\Big(\Delta_n\ge \sqrt{n}/(\log\log n)^2,
 M_n\le \frac {\phi(n)}{\sqrt{1+a-\delta}}\Big)
\le \sum_{j=1}^n \pr\Big(X_j\ne X_{nj}^{\prime}, M_n\le
\frac{\phi(n)}{\sqrt{1+a-\delta}}\Big).
\end{eqnarray*}
and for $1\le j\le n$,
\begin{eqnarray*}
& &\pr\Big(X_j\ne X_{nj}^{\prime}, M_n\le \frac{\phi(n)}{\sqrt{1+a-\delta}}\Big)\\
&\le&\pr\Big(X_j\ne X_{nj}^{\prime},
M_{j-1}\vee \max_{j<k\le n}|S_k-X_j+X_j|\le \frac{\phi(n)}{\sqrt{1+a-\delta}},M_n\le \frac{\phi(n)}{\sqrt{1+a-\delta}}\Big) \\
&\le&\pr\Big(\frac{\sqrt{n}}{\log^p n}<|X_j|
\le 2\frac{\phi(n)}{\sqrt{1+a-\delta}},
M_{j-1}\vee \max_{j<k\le n}|S_k-X_j+X_j|\le \frac{\phi(n)}{\sqrt{1+a-\delta}}\Big) \\
&=&\int_{\frac{\sqrt{n}}{\log^p n}<|y|\le 2\frac{\phi(n)}{\sqrt{1+a-\delta}}}
\pr\Big(M_{j-1}\vee \max_{j<k\le n}|S_k-X_j+y|\le \frac{\phi(n)}{\sqrt{1+a-\delta}}\Big) d \pr(X_j<y).
\end{eqnarray*}
Note that
$M_{j-1}\vee \max_{j<k\le n}|S_k-X_j+y|\overset{\Cal D}=
M_{j-1}\vee \max_{j\le k\le n-1}|S_k+y|$. By Lemma \ref{lem2}, we have
\begin{eqnarray*}
\sup_{|y|\le 2\frac{\phi(n)}{\sqrt{1+a-\delta}}}
 &\pr&\Big(M_{j-1}\vee \max_{j<k\le n}|S_k-X_j+y|
 \le \frac{\phi(n)}{\sqrt{1+a-\delta}}\Big)\\
&\le& C\exp\big\{(-1-a+\delta+\delta^{\prime})\log\log n\big\}.
\end{eqnarray*}
It follows that for $n\in \Cal L$ and $1\le j\le n$,
\begin{eqnarray*}
& & \pr\Big(X_j\ne X_{nj}^{\prime}, M_n\le \frac{\phi(n)}{\sqrt{1+a-\delta}}\Big)\\
&\le& C(\log n)^{-1-a+\delta+\delta^{\prime}}
 \pr\Big(\frac{\sqrt{n}}{\log^p n}<|X_j|
\le 2\frac{\phi(n)}{\sqrt{1+a-\delta}}\Big)\\
&\le& C(\log n)^{-1-a+\delta+\delta^{\prime}}
 \pr\Big(|X|>\sqrt{n}/\log^p n\Big).
\end{eqnarray*}
The remained proof is similar to that of (\ref{eqL3.1}) with Lemma \ref{lem1}(b) instead of Lemma \ref{lem1}(a).

\begin{lemma}\label{lem5}
For any sequence of independent random variables $\{\xi_n; n\ge
1\}$ with mean zero and finite variance, there exists a sequence
of independent normal variables $\{\eta_n; n\ge 1\}$ with $\ep
\eta_n=0$ and $\ep \eta_n^2=\ep \xi_n^2$ such that, for all $Q>2$
and $y>0$,
$$ \pr\Big(\max_{k\le n}|\sum_{i=1}^k \xi_i-\sum_{i=1}^k \eta_i|\ge y\Big)
\le (AQ)^Qy^{-Q}\sum_{i=1}^n \ep |\xi_i|^Q, $$ whenever
$\ep|\xi_i|^Q<\infty$, $i=1,\ldots,n$.
Here, $A$ is a universal
constant.
\end{lemma}

{\bf Proof.} See Sakhaneko (1980,1984, 1985).

\begin{lemma}\label{lem6}
We have that
\begin{eqnarray}\label{eqL6.1}
& &\pr\big(\sup_{0\le s\le 1}|W(s)|\le x-1/(\log\log n)^2\big)-p_n\le\pr\big(M_n^{\ast}\le x\sqrt{B_n}\big)
\non
&\le&
 \pr\big(\sup_{0\le s\le 1}|W(s)|\le x+1/(\log\log n)^2\big)+p_n,
\quad \forall x>0,
\end{eqnarray}
 where $p_n\ge 0$ satisfies
\begin{equation}\label{eqL6.2}
\sum_{n=1}^{\infty}\frac{(\log n)^a(\log\log n)^b}{n} p_n
\le K(a,b,p)<\infty.
\end{equation}
\end{lemma}

{\bf Proof.} By Lemma \ref{lem5}, there exist a universal constant
$A>0$ and a sequence of standard Wiener processes $\{W_n(\cdot)\}$
such that for all $Q>2$,
\begin{eqnarray*}
& & \pr\Big( \max_{k\le n}|S_{nk}^{\ast}-W_n(\frac kn B_n)|
\ge \frac 12 \sqrt{B_n}/(\log\log n)^2\Big) \\
&\le& (AQ)^Q\Big(\frac{(\log\log n)^2}{\sqrt{B_n}}\Big)^Q
\sum_{k=1}^n \ep\big|X_{nk}^{\ast}\big|^Q \\
&\le&  C n \Big(\frac{(\log\log n)^2}{\sqrt{n}}\Big)^Q
 \ep\big[|X|^QI\{|X|\le \sqrt{n}/\log^p n\}\big].
\end{eqnarray*}
On the other hand, by Lemma 1.1.1 of Cs\"org\H o and R\'ev\'esz (1981),
\begin{eqnarray*}
& & \pr\Big(|\max_{0\le s\le 1}|W_n(s B_n)-W_n(\frac{[ns]}{n} B_n)|
\ge \frac 12 \sqrt{B_n}/(\log\log n)^2\Big) \\
&=& \pr\Big(|\max_{0\le s\le 1}|W_n(s)-W_n(\frac{[ns]}{n} )|
\ge \frac 12 \sqrt{\frac 1n}\frac{\sqrt n}{(\log\log n)^2}\Big) \\
&\le& Cn\exp\Big\{-\frac{(\frac 12 \sqrt{n}/(\log\log n)^2 )^2}{3}\Big\}
\le C n\exp\Big\{-\frac 1{12}n/(\log\log n)^4\Big\}.
\end{eqnarray*}
Let
\begin{equation} \label{eqL6.3}
p_n=
\pr\Big(\big|M_n^{\ast}/\sqrt{B_n}-\sup_{0\le s\le 1}|W_n(sB_n)|/\sqrt{B_n}\big|\ge 1/(\log\log n)^2\Big).
\end{equation}
Then $p_n$ satisfies (\ref{eqL6.1}), since $\{W_n(t
B_n)/\sqrt{B_n}; t\ge 0\}\overset{\Cal D}=\{W(t); t\ge 0\}$ for
each $n$. And also,
$$p_n\le C n \Big(\frac{(\log\log n)^2}{\sqrt{n}}\Big)^Q
 \ep\big[|X|^QI\{|X|\le \sqrt{n}/\log^p n\}\big]
+C n\exp\Big\{-\frac 1{12}n/(\log\log n)^4\Big\}. $$
It follows that
\begin{eqnarray*}
& &\sum_{n=1}^{\infty}\frac{(\log n)^a(\log\log n)^b}{n} p_n \\
&\le& K_1+
C\sum_{n=1}^{\infty}\frac{(\log n)^a(\log\log n)^{b+2Q}}{n^{Q/2}}
\ep\big[|X|^QI\{|X|\le \sqrt{n}/\log^p n\}\big] \\
&\le& K_1+
C\sum_{n=1}^{\infty}\frac{(\log\log n)^{b+2Q}}{n(\log n)^{p(Q-2)-a}}
\ep X^2\le K<\infty,
\end{eqnarray*}
whenever $p(Q-2)-a>1$. So, (\ref{eqL6.2}) is satisfied.

%%%%%%%%%%%%%%%%%%%%%%%%%%
\bigskip
Now, we turn to prove Propositions \ref{prop3.1} and \ref{prop3.2}.

{\noindent\bf Proof of Proposition \ref{prop3.1}:}
 Let $0<\delta<\frac 1{2\sqrt{1+a}}$.
Observe that, if $n$ is large enough,
\begin{eqnarray*}
& & \pr\Big\{M_n \le \epsilon \sqrt{\frac{\pi^2 B_n}{8\log\log n}}
 \Big\}\\
&=&\pr\Big\{M_n \le \epsilon
\sqrt{\frac{\pi^2 B_n}{8\log\log n} },
\Delta_n\le \frac{\sqrt{n}}{(\log\log n)^2} \Big\} \\
& & \quad  +\pr\Big\{M_n \le \epsilon \sqrt{\frac{\pi^2
B_n}{8\log\log n} },
\Delta_n> \frac{\sqrt{n}}{(\log\log n)^2} \Big\}\\
&\le&\pr\Big\{M_n^{\ast} \le \epsilon \sqrt{\frac{\pi^2
B_n}{8\log\log n} }
+ \frac{\sqrt{n}}{(\log\log n)^2} \Big\} \\
& & \quad  +\pr\Big\{M_n \le \phi(n)\big(\frac 1{\sqrt{1+a}}+\delta\big),
\Delta_n> \frac{\sqrt{n}}{(\log\log n)^2} \Big\}\\
&\le&\pr\Big\{\sup_{0\le s\le 1}|W(s)| \le \epsilon
\sqrt{\frac{\pi^2 }{8\log\log n} } + \frac{3}{(\log\log n)^2}
\Big\} +p_n +II_n
\end{eqnarray*}
for all $\epsilon\in(1/\sqrt{1+a}-\delta,1/\sqrt{1+a}+\delta)$,
where $II_n$ and $p_n$ are defined in Lemmas \ref{lem4} and
\ref{lem6}, respectively. Similarly, if  $n$ is large enough,
\begin{eqnarray*}
& & \pr\Big\{M_n \le \epsilon \sqrt{\frac{\pi^2 B_n}{8\log\log n}
}
\Big\}\\
&\ge&\pr\Big\{M_n \le \epsilon
\sqrt{\frac{\pi^2 B_n}{8\log\log n}
}, \Delta_n\le \frac{\sqrt{n}}{(\log\log n)^2} \Big\} \\
&\ge&\pr\Big\{M_n^{\ast} \le \epsilon \sqrt{\frac{\pi^2
B_n}{8\log\log n} }- \frac{\sqrt{n}}{(\log\log n)^2},
\Delta_n\le \frac{\sqrt{n}}{(\log\log n)^2}  \Big\} \\
&\ge&\pr\Big\{M_n^{\ast} \le \sqrt{B_n} \big[\epsilon
\sqrt{\frac{\pi^2 }{8\log\log n} } - \frac{2}{(\log\log n)^2}\big]
\Big\}
\\
& &\quad -\pr\Big\{M_n^{\ast} \le \sqrt{B_n} \big[\epsilon
\sqrt{\frac{\pi^2 }{8\log\log n} } - \frac{2}{(\log\log
n)^2}\big],
 \Delta_n> \frac{\sqrt{n}}{(\log\log n)^2}\Big\} \\
&\ge&\pr\Big\{\sup_{0\le s\le 1}|W(s)| \le \epsilon
\sqrt{\frac{\pi^2 }{8\log\log n} }
 - \frac{3}{(\log\log n)^2} \Big\} \\
& & \quad  -\pr\Big\{M_n^{\ast} \le \phi(n)\big(\frac 1{\sqrt{1+a}}+\delta\big),
\Delta_n> \frac{\sqrt{n}}{(\log\log n)^2} \Big\}\\
&\ge&\pr\Big\{\sup_{0\le s\le 1}|W(s)| \le \epsilon
\sqrt{\frac{\pi^2 }{8\log\log n} } - \frac{3}{(\log\log n)^2}
\Big\} -p_n -I_n
\end{eqnarray*}
for all $\epsilon\in(1/\sqrt{1+a}-\delta,1/\sqrt{1+a}+\delta)$,
where $I_n$ is  defined in Lemma \ref{lem3}. Choosing $\delta>0$
small enough and letting $q_n=p_n+I_n+II_n$ complete the proof by
Lemmas \ref{lem3}, \ref{lem4} and \ref{lem6}.

\bigskip

%%%%%%%%%%%%%%%%%%%%%%%%%

{\noindent\bf Proof of Proposition \ref{prop3.2}:} Let
$\{W_n(\cdot)\}$ be a sequence of standard Wiener processes being
defined in the proof of Lemma \ref{lem6}, and let $p_n$ be defined
in (\ref{eqL6.3}). And set
$$q_n^{\ast}=
\pr\Big(\big|M_n/\sqrt{B_n}-\sup_{0\le s\le 1}|W_n(sB_n)|/\sqrt{B_n}\big|\ge 3/(\log\log n)^2\Big). $$
Then $q_n^{\ast}$ satisfies (\ref{eqprop3.2.1}), and
$$q_n^{\ast}\le \pr\big(\Delta_n\ge \sqrt{n}/(\log\log n)^2\big)+p_n. $$
By Lemma \ref{lem6},
$$\sum_{n=1}^{\infty}
\frac{(\log\log n)^b}{n\log n} p_n \le K_1(b,p)<\infty.
$$
Also, following the lines in the proof of (\ref{eqL3.1}) we have
\begin{eqnarray*}
& &\sum_{n=1}^{\infty}\frac{(\log\log n)^b}{n\log n}
   \pr\big(\Delta_n\ge \sqrt{n}/(\log\log n)^2\big) \\
&\le&\sum_{n\in \Cal L}\frac{(\log\log n)^b}{n\log n}\cdot n
   \pr\big(|X|> \sqrt{n}/(\log n)^p\big)
+\sum_{n\not\in \Cal L}\frac{(\log\log n)^{b+2}}{n^{3/2}\log n}
  \beta_n\\
&\le&\sum_{n=1}^{\infty}\frac{(\log\log n)^b}{\log n}
   \pr\big(|X|> \sqrt{n}/(\log n)^p\big) \\
& & \quad
+\sum_{n=1}^{\infty}\frac{(\log\log n)^{b+2}}{\sqrt{n}\log n}
   \ep\big[|X|I\{|X|> \sqrt{n}/(\log n)^p\}\big]
\\
&\le&C\ep\big[X^2(\log|X|)^{2p-1}(\log\log |X|)^b\big]
  +C\ep\big[X^2(\log|X|)^{p-1}(\log\log |X|)^{b+2}\big]\\
&\le& C\ep X^2<\infty.
\end{eqnarray*}
So, $q_n^{\ast}$ satisfies (\ref{eqprop3.2.2}).

%%%%%%%%%%%%%%%%%%%%%%%%%%%%%%%%%%%%%%%%%%%%%%
%%%%%%%%%%%%%%%%%%%%%%%%%%%%%%%%%%%%%%%%%%%%

\section{Proofs of the Theorems.}
\setcounter{equation}{0}

%%%%%%%%%%%%%%%
\subsection{Proofs of the direct parts.}

Without losing of generality, we assume that $\ep X=0$ and $\ep X^2=1$.

{\noindent\bf Proof of the direct part of Theorem \ref{th1}:}
Let $\delta>0$ small enough and $\{q_n\}$ be such that (\ref{eqprop3.1.1}) and (\ref{eqprop3.1.2}) hold. Then
\begin{eqnarray*}
 \lim_{\epsilon\nearrow 1/\sqrt{1+a}}
(\frac{1}{\sqrt{1+a}}-\epsilon)^{b+1} \sum_{n=1}^{\infty}
\frac{(\log n)^a (\log\log n)^b}{n} q_n=0,
\end{eqnarray*}
by (\ref{eqprop3.1.2}). Notice that $a_n(\epsilon)\to 0$.
By (\ref{eqprop3.1.1}), we have that for $n$ large enough,
\begin{eqnarray*}
&&\pr\Big\{\sup_{0\le s\le 1}|W(s)| \le
\sqrt{\frac{\pi^2 }{8\log\log n} }
(\epsilon+a_n(\epsilon)) - \frac{3}{(\log\log n)^2} \Big\}
-q_n
\\
&\le& \pr\Big\{M_n \le
\sqrt{\frac{\pi^2 B_n}{8\log\log n} }
(\epsilon+a_n(\epsilon)) \Big\}
\\
&\le&\pr\Big\{\sup_{0\le s\le 1}|W(s)| \le
\sqrt{\frac{\pi^2 }{8\log\log n} }
(\epsilon+a_n(\epsilon)) + \frac{3}{(\log\log n)^2} \Big\}
+q_n,
\\
& &  \qquad \forall \epsilon\in
(\frac 1{\sqrt{1+a}}-\delta/2, \frac 1{\sqrt{1+a}}+\delta/2).
\end{eqnarray*}
On the other hand, by Proposition \ref{prop2.1},
\begin{eqnarray*}
& & \lim_{\epsilon\nearrow 1/\sqrt{1+a}}
(\frac{1}{\sqrt{1+a}}-\epsilon)^{b+1}
\sum_{n=1}^{\infty} \frac{(\log n)^a (\log\log n)^b}{n} \\
& & \qquad \quad \cdot \pr\Big\{\sup_{0\le s\le 1}|W(s)| \le
\sqrt{\frac{\pi^2 }{8\log\log n} }
(\epsilon +a_n(\epsilon))\pm \frac{3}{(\log\log n)^2} \Big\}
\\
& & \quad=\frac 4{\pi}(\frac{1}{2(1+a)^{3/2}})^{b+1}\Gamma(b+1)
\exp\Big\{2(1+a)^{3/2}\tau\Big\}.
\end{eqnarray*}
It follows that
\begin{eqnarray}\label{eqP1.1}
& & \lim_{\epsilon\nearrow 1/\sqrt{1+a}}
(\frac{1}{\sqrt{1+a}}-\epsilon)^{b+1} \sum_{n=1}^{\infty}
\frac{(\log n)^a (\log\log n)^b}{n} \pr\Big\{M_n \le
\sqrt{\frac{\pi^2 B_n}{8\log\log n} } (\epsilon
+a_n(\epsilon))\Big\}
\non
& & \qquad \qquad \qquad
=\frac4{\pi}(\frac{1}{2(1+a)^{3/2}})^{b+1}\Gamma(b+1)
 \exp\Big\{2(1+a)^{3/2}\tau\Big\}.
\end{eqnarray}

Finally, noticing the condition (\ref{co1}), we have
$$0\le n-B_n\le 2 n \ep[X^2I\{|X|\ge \sqrt{n}/\log^p n\}]
=o(n(\log\log n)^{-1}). $$
Let $a_n^{\prime}(\epsilon)=\sqrt{n/B_n}(\epsilon+a_n(\epsilon))-\epsilon$.
Then
$$\pr\Big\{M_n \le \phi(n)(\epsilon +a_n(\epsilon))\Big\}
=\pr\Big\{M_n \le
\sqrt{\frac{\pi^2 B_n}{8\log\log n} }
(\epsilon +a_n^{\prime}(\epsilon))\Big\}, $$
and,
$$a_n^{\prime}(\epsilon)\log\log n
=\epsilon \frac{(n-B_n)\log\log n}{\sqrt{B_n}(\sqrt n+\sqrt{B_n})}+\sqrt{\frac{n}{B_n}}a_n(\epsilon)\log\log n
\to \tau $$
as $n\to \infty$ and $\epsilon\nearrow 1/\sqrt{1+a}$. Now, (\ref{eqT1.1}) follows from (\ref{eqP1.1}).

%%%%%%%%%%%%%%%

\bigskip

{\noindent \bf Proof of the direct part of Theorem \ref{th2}:}
Noticing $B_n\sim n$ and Proposition \ref{prop3.2}, for any
$0<\delta<1$ we have for $n$ large enough and all $\epsilon>10^3$,
\begin{eqnarray*}
& & \pr\Big\{\sup_{0\le s\le 1}|W(s)|\le (1-\delta)\epsilon
\sqrt{\frac{\pi^2}{8\log\log n}}\Big\}-q_n^{\ast} \\
&\le& \pr\Big\{\sup_{0\le s\le 1}|W(s)|\le \epsilon
\sqrt{\frac{\pi^2}{8\log\log n}} - 3/(\log\log n)^2\Big\}
  -q_n^{\ast} \\
&\le&\pr\Big\{M_n\le \epsilon \sqrt{\frac{\pi^2 B_n }{8\log\log n}} \Big\}\\
&\le&\pr\Big\{M_n\le \phi(n) \epsilon\Big\} \le \pr\Big\{M_n\le
(1+\frac{\delta}2) \epsilon
 \sqrt{\frac{\pi^2 B_n }{8\log\log n}}
\Big\} \\
&\le&\pr\Big\{\sup_{0\le s\le 1}|W(s)|\le (1+\frac{\delta}2)
\epsilon \sqrt{\frac{\pi^2}{8\log\log n}}
  + 3/(\log\log n)^2\Big\}
  +q_n^{\ast} \\
&\le&\pr\Big\{\sup_{0\le s\le 1}|W(s)|\le (1+\delta) \epsilon
\sqrt{\frac{\pi^2}{8\log\log n}} \Big\}
  +q_n^{\ast}.
\end{eqnarray*}
So, by Propositions \ref{prop2.2} and \ref{prop3.2},
\begin{eqnarray*}
& &(1-\delta)^{2(b+1)}\frac4{\pi}\Gamma(b+1)
\sum_{k=0}^{\infty}\frac{(-1)^k}{(2k+1)^{2b+3}}\\
&=& \lim_{\epsilon\nearrow \infty}
\epsilon^{-2(b+1)} \sum_{n=1}^{\infty}
\frac{(\log\log n)^b}{n\log n}
 \pr\Big\{\sup_{0\le s\le 1}|W(s)|\le (1-\delta)\epsilon
\sqrt{\frac{\pi^2}{8\log\log n}}\Big\} \\
&\le& \liminf_{\epsilon\nearrow \infty}
\epsilon^{-2(b+1)} \sum_{n=1}^{\infty}
\frac{(\log\log n)^b}{n\log n}
\pr\Big\{M_n\le \phi(n) \epsilon\Big\} \\
&\le& \limsup_{\epsilon\nearrow \infty}
\epsilon^{-2(b+1)} \sum_{n=1}^{\infty}
\frac{(\log\log n)^b}{n\log n}
\pr\Big\{M_n\le \phi(n) \epsilon\Big\} \\
&\le& \lim_{\epsilon\nearrow \infty}
\epsilon^{-2(b+1)} \sum_{n=1}^{\infty}
\frac{(\log\log n)^b}{n\log n}
 \pr\Big\{\sup_{0\le s\le 1}|W(s)|\le (1+\delta)\epsilon
\sqrt{\frac{\pi^2}{8\log\log n}}\Big\} \\
&= &(1+\delta)^{2(b+1)}\frac4{\pi}\Gamma(b+1)
\sum_{k=0}^{\infty}\frac{(-1)^k}{(2k+1)^{2b+3}}.
\end{eqnarray*}
Letting $\delta\to 0$ completes the proof.

%%%%%%%%%%%%%%%
\bigskip
{\noindent\bf Proof of the direct part of Theorem \ref{th3}:}
 Choose $\delta>0$ small enough and $\{q_n\}$
 for (\ref{eqprop3.1.1}) and (\ref{eqprop3.1.2}) to hold.
 By a standard argument (see Feller 1945), we can assume that
$$ \frac{1}{\sqrt{\log\log n}}\big(\frac{1}{\sqrt{1+a}}-\frac{\delta}2\big)
\le \psi(n)\le
\frac{1}{\sqrt{\log\log n}}\big(\frac{1}{\sqrt{1+a}}+\frac{\delta}2\big). $$
Then for $n$ large enough,
\begin{eqnarray*}
& &\pr\Big\{\sup_{0\le s\le 1}|W(s)| \le
\sqrt{\pi^2/8 }\big[\psi(n)-5/(\log\log n)^2\big]\Big\}
-q_n \\
&\le&\pr\Big\{\sup_{0\le s\le 1}|W(s)| \le
\sqrt{\pi^2/8 }\psi(n)-3/(\log\log n)^2\Big\}
-q_n \\
&\le& \pr\Big\{M_n \le
\sqrt{\pi^2 B_n /8 }\psi(n)\Big\}\\
&\le& \pr\Big\{\sup_{0\le s\le 1}|W(s)| \le
\sqrt{\pi^2/8 }\psi(n)+3/(\log\log n)^2\Big\}
+q_n \\
&\le& \pr\Big\{\sup_{0\le s\le 1}|W(s)| \le
\sqrt{\pi^2/8 }\big[\psi(n)+5/(\log\log n)^2\big]\Big\}
+q_n.
\end{eqnarray*}
Notice that $J_{ab}\big(\psi(n)\pm 5/(\log\log n)^2\big)<\infty$ if and only if $J_{ab}(\psi(n))<\infty$.
So, by Propositions \ref{prop2.3} and \ref{prop3.1} we have
\begin{eqnarray}\label{eqP3.1}
& & \sum_{n=1}^{\infty} \frac{(\log n)^a (\log\log
n)^b}{n}\pr\big\{M_n \le \sqrt{\pi^2 B_n /8 }\psi(n)\big\}<\infty
\text{ or } =\infty \non
 & & \text{ according as }
J_{ab}(\psi)<\infty \text{ or } = \infty.
\end{eqnarray}

Finally, if condition (\ref{co2}) is satisfied, then
$$0\le n-B_n=O(n/\log\log n), $$
which implies that
$$\sqrt{n}\psi(n)=\sqrt{B_n}\psi(n)\sqrt{1+O(1/\log\log n)}. $$
Notice that
$J_{ab}\big(\psi(n)\sqrt{1\pm O(1/\log\log n)}\big)<\infty$
if and only if $J_{ab}(\psi(n))<\infty$. (\ref{eqT3.1}) follows from (\ref{eqP3.1}).

%%%%%%%%%%%%%%%%%%%%
%%%%%%%%%%%

\subsection{Proofs of the converse parts.}
Now, we turn to prove the converse parts of Theorem
\ref{th1}-\ref{th3}. The proof is organized as follows. First, we
show that $\ep X^2=\infty$ is impossible. Secondly, we show that
$\ep X\ne 0$ is impossible if $\ep X^2<\infty$. Thirdly, we show
that $\ep X^2\ne \sigma^2$ is impossible if $\ep X=0$ and $\ep
X^2<\infty$. At last, we show (\ref{eqT1.2}) and (\ref{eqT3.2}).

From Esseen (1968) (see also Petrov 1995) it is easy to see that
for all $m\ge 1$,
$$\pr(|S_m|\le 2\sqrt{m})\le K\big(\int_{-2\sqrt{m}}^{2\sqrt{m}}d \widetilde{F}(x)\big)^{-1/2}, $$
where $\widetilde{F}(x)$ is the distribution function of the
symmetrized $X$, and $K$ is an absolute constant.
So, if $\ep
X^2=\infty$, then for any $M>2$ we can choose $m_0\ge 9$ large
enough such that
$$\pr(|S_m|\le 2\sqrt{m})\le e^{-2M}, \quad m\ge m_0. $$
For $\epsilon>0$, we let $m=[\epsilon^2n/\log\log n]$, and
$N=[n/m]$, then for all $n\ge m_0^2$ and $\epsilon\ge 1$,
\begin{eqnarray}\label{eq3.1}
\pr\big(M_n\le \epsilon (n/\log\log n)^{1/2}\big)
&\le&\pr\big(|S_{km}-S_{(k-1)m}|\le 2\sqrt{m}, k=1,\ldots, N\big)
\non
&\le& e^{-2MN}\le \exp\big\{-M\frac{\log\log n}{\epsilon^2}\big\}.
\end{eqnarray}
By this inequality, for any $a$, $b$ and $0<\epsilon_1<\epsilon_2<\infty$
there exists a constant $C=C(a,b,\epsilon_1,\epsilon_2)$ for which
$$ \sup_{\epsilon\in (\epsilon_1,\epsilon_2)}
\sum_{n=1}^{\infty}\frac{(\log n)^a(\log\log n)^b}{n}
 \pr\big(M_n\le \epsilon (n/\log\log n)^{1/2}\big)\le C<\infty,$$
which implies that (\ref{eqT1.1}) and (\ref{eqT3.1}) can not hold.
Also, by (\ref{eq3.1}) and (\ref{eq2.3}), for any $b>-1$ we have
\begin{eqnarray*}
& &\limsup_{\epsilon\nearrow \infty}
 \epsilon^{-2(b+1)}
\sum_{n=1}^{\infty}\frac{(\log\log n)^b}{n\log n}
 \pr\big(M_n\le \epsilon (n/\log\log n)^{1/2}\big)\\
& & \le \lim_{\epsilon\nearrow \infty}
 \epsilon^{-2(b+1)}
\sum_{n=1}^{\infty}\frac{(\log\log n)^b}{n\log n}
\exp\big\{-M\frac{\log\log n}{\epsilon^2}\big\} \\
& & =M^{-(b+1)}\Gamma(b+1)\to 0 \quad \text{as}\; M\to \infty,
\end{eqnarray*}
which implies that (\ref{eqT2.1}) can not hold.

If $\ep X^2<\infty$ and $\ep X=\mu\ne 0$, then
$$\pr(M_n\le \frac{n|\mu|}{2})\le \pr(|S_n|\le \frac{n|\mu|}{2})
\le \pr(|S_n-n\mu|\ge  \frac{n|\mu|}{2})\le \frac{4\ep X^2}{\mu^2 n^2}. $$
It follows that
\begin{eqnarray*}
\sum_{n=1}^{\infty}\frac{(\log n)^a(\log\log n)^b}{n}\pr(M_n\le \frac{n|\mu|}{2})
 \le \frac{4\ep X^2}{\mu^2}\sum_{n=1}^{\infty}\frac 1{n^3}
\le \frac{8\ep X^2}{\mu^2}.
\end{eqnarray*}
It follows that for any $\epsilon>0$,
\begin{eqnarray*}
& &\sum_{n=1}^{\infty}\frac{(\log n)^a(\log\log n)^b}{n}
\pr(M_n\le \epsilon(n/\log\log n)^{1/2}) \\
&\le&\sum_{n<(2\epsilon /|\mu|)^2}\frac{(\log n)^a(\log\log n)^b}{n}
\pr(M_n\le \epsilon(n/\log\log n)^{1/2}) \\
& &+\sum_{n\ge (2\epsilon /|\mu|)^2 }
 \frac{(\log n)^a(\log\log n)^b}{n}
\pr(M_n\le \epsilon(n/\log\log n)^{1/2}) \\
&\le&\sum_{n\ge (2\epsilon /|\mu|)^2}
\frac{(\log n)^a(\log\log n)^b }{n}
\pr(M_n\le \frac{n|\mu|}{2})
+\sum_{n< (2\epsilon /|\mu|)^2}\frac{(\log n)^{|a|+|b|}}{n}
 \\
&\le&
\frac{8\ep X^2}{\mu^2}
+4\Big(\log(\frac{2\epsilon}{|\mu|})^4\Big)^{|a|+|b|+1},
\end{eqnarray*}
which implies none of (\ref{eqT1.1}), (\ref{eqT2.1}) and
(\ref{eqT3.1}) can hold.

Suppose that $\ep X=0$ and $\ep X^2<\infty$, and that
(\ref{eqT1.1}), (\ref{eqT2.1}) and (\ref{eqT3.1}) hold for some
positive constant $\sigma$. By the direct part of Theorem
\ref{th2}, (\ref{eqT2.1}) shall hold with $\ep X^2$ taking the
place of $\sigma^2$. This is obvious a contradiction if $\ep
X^2\ne \sigma^2$. Notice that (\ref{eqP1.1}) and (\ref{eq3.1})
hold whenever $\ep X=0$ and $\ep X^2<\infty$. However, if $\ep
X^2\ne \sigma^2$, then  (\ref{eqP1.1}) and (\ref{eqP3.1}) are
contradictory to (\ref{eqT1.1}) and (\ref{eqT3.1}) respectively,
since $B_n\sim n \ep X^2$.

Finally, we show (\ref{eqT1.2}) and (\ref{eqT3.2}). Suppose that
(\ref{eqT1.2}) fails. Without losing of generality, we can assume
that $\sigma^{-2}\ep [X^2I\{|X|\ge \sqrt{n}/(\log n)^p\}]\ge
\tau_0/\log\log n$ for some $\tau_0>0$ and all $n\ge 1$. Then
$n\sigma^2-B_n\ge n\ep [X^2I\{|X|\ge \sqrt{n}/(\log n)^p\}]\ge
n\sigma^2\tau_0/\log\log n$. Let
$a_n^{\prime}(\epsilon)=\sqrt{1+\tau_0/\log\log
n}\big(\epsilon+a_n(\epsilon)\big)-\epsilon$. Then
$$ a_n^{\prime}(\epsilon)\log\log n\to \tau+\tau_0/(2\sqrt{1+a}),$$
and
$$\pr\Big\{M_n\le \sigma\phi(n)\big(\epsilon+a_n(\epsilon)\big)\Big\}
 \ge \pr\Big\{M_n\le \sqrt{\frac{\pi^2 B_n}{8\log\log n}}
 \big(\epsilon+a_n^{\prime}(\epsilon)\big)\Big\}. $$
It follows that (\ref{eqT1.1}) is contradictory to (\ref{eqP1.1}).

Now, suppose that (\ref{eqT1.2}) fails. Let $d_n=(\log\log
n)\ep[X^2I\{|X|\ge \sqrt{n}/(\log n)^p\}]/\sigma^2$. Then
$d_n/\log\log n$ is non-increasing in $n$ and $d_n\to \infty$. So,
one can find a non-decreasing sequence $\{b_n\}$ of positive
numbers for which
 $e\le b_n\to \infty$,
$$ \sum_{n=1}^{\infty}\frac{1}{n(\log n)(\log\log n)b_n}<\infty \;
\text{ and } \;
\sum_{n=1}^{\infty}\frac{d_n}{n(\log n)(\log\log n)b_n}=\infty. $$
Define
$$\psi(n)=1/\sqrt{(1+a)\log\log n+(1+b)\log\log\log n+\log b_n}$$
and
$$\widetilde{\psi}(n)= \psi(n)/\sqrt{1-d_n/(2\log\log n)}. $$
Then
$$J_{ab}(\psi)=\sum_{n=1}^{\infty}\frac{1}{n(\log n)(\log\log n)b_n}<\infty$$
and
\begin{eqnarray*}
J_{ab}(\widetilde{\psi})
&=& \sum_{n=1}^{\infty}
\frac{(\log n)^a(\log\log n)^b}{n}\exp\Big\{-\frac{1-d_n/(2\log\log n)}{\psi^2(n)} \Big\}
\\
&\ge& \sum_{n=1}^{\infty}
\frac{(\log n)^a(\log\log n)^b}{n}\exp\Big\{-\frac{1}{\psi^2(n)} \Big\}
\exp\Big\{\frac{1+a}{2}d_n\Big\} \\
&\ge& c \sum_{n=1}^{\infty}\frac{d_n}{n(\log n)(\log\log n)b_n}=\infty.
\end{eqnarray*}
However,
$$\pr \Big\{ M_n\le \sqrt{ \frac{\pi^2\sigma^2n }{8} }\psi(n) \Big\}
=\pr\Big\{ M_n\le \sqrt{ \frac{\pi^2B_n}{8} } \psi(n) \sqrt{\frac{n\sigma^2}{B_n}}\Big\}
\ge \pr\Big\{ M_n\le \sqrt{ \frac{\pi^2B_n}{8} }
 \widetilde{\psi}(n) \Big\},
$$
since $n\sigma^2-B_n\ge n\sigma^2 d_n$.
It follows that (\ref{eqT3.1}) is contradictory to (\ref{eqP3.1}).
The proof is now completed.
%%%%%%%%%%%%%%%%
\newpage

%%%%%%%%%%%%%%%%%%%%%%%%%%%%%%%%%%%%%%%%%%%%%
%\begin{center}{\bf References}\end{center}

\end{document}